\documentclass[12pt]{article}
\usepackage{amsmath, amssymb, amsfonts}
\usepackage{geometry}                		
\usepackage{graphicx}				
\usepackage{amssymb}
 at 10truept
\font \symb = msbm10 at 12truept

 at 16truept
 at 16truept
 at 12truept

 at 10truept
 at 10truept
 at 8truept
 at 10truept
 at 8truept
 at 10 truept
 at 11 truept
 at 9truept
 at 12truept
 at 12truept
 at 7truept
 at 14 truept
 at 12truept
 at 9truept
 at 14truept
 at 17 truept
\newcommand\proof{\smallbreak\noindent {\sc Proof $-\;$}}

\newcommand\R{\hbox {\symb R}}

\def\XXint#1#2#3{{\setbox0=\hbox{$#1{#2#3}{\int}$ }
		\vcenter{\hbox{$#2#3$ }}\kern-.6\wd0}}

\def\abstr{\if@twocolumn
	\section*{abstr}
	\else \small
	\quotation\fi
	{\bf Abstract.}
}
\def\endabstr{\if@twocolumn\else\endquotation\fi}
\mark{{}{}}
\begin{document}
	
	\title{{\LARGE {On Existence of  Two Positive Solutions for  Nonlinear Subelliptic Equations Involving Nonuniformly P-Laplacian
		}}\\
		[30pt]}
	\author{ Farman Mamedov$^{*}$ 
	and Jasarat Gasimov$^{**}$ \\
	{\small {(mfamannn@gmail.com; jasarat.gasimov@emu.edu.tr;  
			)
	}}\\
	[1cm] $^*$Mathematics and Mechanics Institute \\
	of National Academy of Sciences. \\
	(\emph{\ B.Vahabzade str., 9, AZ 1141, Baku, Azerbaijan}) \\
	$^{**}$Department of Mathematics, Eastern Mediterranean University,\\
	Mersin 10, 99628, T.R. North Cyprus, Turkey.
	\\
	\\
}
\maketitle

\newtheorem{teo}{\quad Theorem}
\newtheorem{prop}[teo]{Proposition} %
\newtheorem{defi}[teo]{Definition}
\newtheorem{lem}[teo]{\quad Lemma} %
\newtheorem{cor}[teo]{Corollary}
\newtheorem{rem}[teo]{\quad Remark} %
\newtheorem{ex}[teo]{Example}
\newtheorem{cla}{\quad Claim}

\pagestyle{myheadings} \baselineskip 13.5pt

\bigskip \bigskip

\thispagestyle{empty} \medskip \indent

\begin{abstr} \small{ In this paper, we study an existence result of two different positive solutions for the problem 
		$$
		\mbox {div } \left ( \vert \nabla_\omega u\vert^{p-2}\nabla_\omega u \right )+v(x) u^{q-1}+\mu u^{\gamma-1}=0, \quad z\in  \Omega, \quad u \Big \vert_{\partial \Omega}=0.
		$$ 
		where $\nabla_\omega u =\left ( \omega^{1/p} \nabla_x, \, \nabla_y  \right ) $ and 
		$
		\vert \nabla_\omega u\vert=  \left (\omega(x)^{2/p} \vert  \nabla_{x}u \vert ^2+\vert  \nabla_{y}u\vert^2 \right )^{\frac{1}{2}}
		$ 
		for all  $z \in \Omega,  \  $ with $ z =(x, y ), \, x \in \R^n, \, y \in \R^m$ and $\Omega\subset \R^N$ is a bounded domain, $ N=n+m, \, p\in(1, N).$
		We restrict on the range of nonlinearities $ q \in (p, pN/(N-p)) $ and $ \gamma \in \left (1, N/(N-1)\right ) $ (or $\gamma\in (1, p)$ and $v^{-\gamma/(q-\gamma)}\in L_{1,loc}(\Omega)$ additionally) and $ \mu \in (0, \Lambda) $ for sufficiently small $ \Lambda; $ the $x$ variable dependent weight functions $ v \in A_\infty, \, \omega \in A_p $  are of the
		Muckenhoupt classes and satisfy also to the balance condition of Chanillo-Wheeden's type.}
\end{abstr}

\bigskip

\bigskip

\bigskip


\medskip

\vspace{2cm} \quad

\bigskip

\noindent

\section{Introduction}

\qquad

The focus of this work is on the Dirichlet problem for the class of second-order nonlinear equations where the major portion is a nonuniformly $p$-Laplacian operator of N variables and the subordinate component has nonlinearities of the weighted convex and concave function types with the parameter
\begin{equation}\label{eq1}
	\mbox{div} \left ( \vert  \nabla_\omega u  \vert ^{p-2} \nabla _\omega u \right )+v(x) u^{q-1} + \mu u^{\gamma-1}=0  \quad z\in  \Omega
\end{equation} 
where $z=(x,y)\in \Omega$\quad  and \quad $\Omega\subset R^{N}, N=n+m$ is bounded domain.
$\nabla_\omega u =\left ( \omega^{1/p} \nabla_x, \, \nabla_y  \right ) $ is the non-uniformly degenerate gradient, $ \nabla_x=\left (\partial / \partial x_1,\cdots \partial / \partial x_n \right )$ and  $ \nabla_y=\left (\partial / \partial y_1,\cdots \partial / \partial y_m \right ). $ 

Assume  $ p\in (1, N), \,   q \in (p, pN/(N-p)) $ and $ \gamma \in \left (1, N/(N-1)\right ) $ (or $\gamma\in (1, p)$ and $v^{-\gamma/(q-\gamma)}\in L_{1,loc}(\Omega)$ additionally),  $ \mu \in (0, \Lambda) $ for sufficiently small $ \Lambda, $ the weight functions $ v(x) \in A_\infty, \, \omega(x) \in A_p $  are of proper Muckenhoupt classes of $x$-variable, also satisfying the Chanillo-Wheeden's type balance condition. The principal component of equation \eqref{eq1} does not constitute a uniformly elliptic operator, mainly due to the fact that the function $\omega(x)$ from the $A_p $ class is, in general, unbounded both above and below. Hence, we characterize equation \eqref{eq1} as nonuniformly $p$-Laplacian.

In a more detailed manner, this paper is dedicated to investigating the Dirichlet problem within the class of nonuniformly $p$-Laplacian equations. These equations involve a nonlinearity expressed as the combination of nonlinearities originating from a weighted convex function and a concave function with a parameter.
\begin{equation}\label{1c}
	\begin{split}
		& \mbox{div} \left ( \vert  \nabla_\omega u  \vert ^{p-2} \nabla _\omega u \right )+v(x) u^{q-1} + \mu u^{\gamma-1}=0\quad \text{in} \quad \Omega, \\
		& u>0 \quad \text{in} \quad \Omega, \\
		& u=0 \quad \text{in} \quad \partial \Omega .
	\end{split}
\end{equation}
A function that satisfies problem (4) in a weak sense is what we refer to as a solution $u \in  \mathring { \mathcal W}_{p}^1 \left ( \Omega ,\omega (x)dz \right )$  satisfying the integral identity $ \forall \phi \in \mathring { \mathcal W}_{p}^1 \left ( \Omega ,\omega (x)dz \right ): $
\begin{equation}\label{2}
	\int\limits_{\Omega} 
	\vert \nabla_\omega u\vert^{p-2}\nabla_\omega u\nabla_\omega \phi \, dz
	- \int\limits_\Omega v(x) u _+^{q-1} \phi \, dz- \mu\int\limits_\Omega  u _+^{\gamma-1} \phi \, dz=0 .
\end{equation}
(for the definition of the space $\mathring { \mathcal W}_{p}^1 \left ( \Omega ,\omega (x)dz \right )$), see\cite{FC}).

The search for positive solutions to the Dirichlet problem has been a focal point of recent research, employing topological methods, notably the mountain-pass theorem, to address nonlinear problems. In this context, the paper \cite{ABC} merits attention, as it delves into the case involving nonlinearity expressed as the sum of convex and concave functions. The paper specifically tackles the following problem ( case $p=2$):
\begin{equation}\label{0}
	\begin{split}
		& \Delta_p u+\lambda u^{q-1}+ u^{\gamma-1}=0\quad \text{in} \quad \Omega, \\
		& u>0 \quad \text{in} \quad \Omega, \\
		& u=0 \quad \text{on} \quad \partial \Omega,
	\end{split}
\end{equation} 
where $\Omega\subset \R^N$ is bounded ,  $N \geq 3, \lambda >0, \, 1<\gamma<p<q \leq pN/(N-p), \Delta_p u $- denotes the p-Laplacian, and $\Delta_p u:=\mbox{div} \left ( \vert\nabla u\vert^{p-2}\nabla u\right ) $ for all $u\in \mathring { \mathcal W}_{p}^1 \left ( \Omega  \right )$ was considered on case $p=2$ by  Ambrosetti, Brezis and Cerami.They demonstrated that for any  $\lambda \in (0,\lambda_{0})$ for some positive number $\lambda_{0}>0$, there are at least two distinct positive solutions to this issue in \cite{ABC} ( see also  \cite{LD} ).In \cite{Alo} and \cite{QZ}, the study \cite{ABC} was extended to nonlinear cases with p-Laplacian in the core component.In general cases, the problem were studied in \cite{Catr},\cite{FR},\cite{Hu}, \cite{MoP}, \cite{PoR}. In \cite{aa}  the problem considered
\begin{equation}\label{AA}
	\begin{split}
		& -\Delta_p u=\lambda \vert u\vert^{q-2}u+\vert u\vert^{p^{*}-2}u\quad \text{in} \quad \Omega, \\
		& u>0 \quad \text{in} \quad \Omega, \\
		& u=0 \quad \text{on} \quad \partial \Omega
	\end{split}
\end{equation}
with $\Omega\subset \R^N$ is a smooth bounded domain, $1<q<p<N,  \lambda>0,  p^{*}=Np/(N-p). $ It was proved that there exists a $\lambda_{0}$ such that for all $0<\lambda<\lambda_{0}$, this problem has at least two positive solutions, where either $2N/(N+2)<p<3$ and $1<q<p$, or $p\geq 3$ and $p>q>p^{*}-2/(p-1).$ The problem \eqref{eq1} is studied also in \cite{FC}, on the case of  the principal part being a linear nonuniformly elliptic operator.

This study attempts to achieve an appropriate result for the class of equations of type \eqref{eq1} that are nonuniformly elliptic. We need new nonuniform inequalities of Poincare$-$Sobolev (PS) type (\cite{FM5}) and their compact counterparts in order to obtain a solution of problem \eqref{1c} based on the use of the mountain-pass Theorem [13] (cf.\cite{YJ}).

We have all the tools we need to use the mountain-pass theorem and variational theory to prove that there is a solution to problem \eqref{1c}, provided that the proper Muckenhoupt conditions $v \in A_{\infty}, \omega\in A_{p}$ are assumed on the degeneration in the principal part $\omega (x)$ and on the weight $v(x)$ in the nonlinear term $v(x)u^{q-1},$ the balance condition \eqref{BC}, and the compactness condition \eqref{CC} below.We restrict on the range of nonlinearities $ q \in (p, pN/(N-p)) $ and $ \gamma \in \left (1, N/(N-1)\right ) $ (or $\gamma\in (1, p)$ and $v^{-\gamma/(q-\gamma)}\in L_{1,loc}(\Omega)$ additionally).

\smallskip
\section{Preliminaries}\label{s2}

\quad

For  $ 1<p<\infty ,$ the $ p^\prime $ is denoted conjugate number, $ \frac{1}{p}+\frac{1}{p^\prime}=1 $; for $  p=\infty ,$ we put $ p^\prime=1 $, and for $ p^\prime=\infty ,$  $ p=1 $.Let $ \chi_{E} $-denote the characteristic function of set $ E $.
$ C, C_{1}, C_{2}, ... $ are denoted positive different constants which can vary and whose values are not essential(see e.g. \cite{FC}).

\smallskip

Recall some basic notions relate to Sobolev spaces defined above (see, e.g. \cite {VM}). Let $\omega: \R^n\to (0,\infty)$ be a positive measureable function satisfying the $A_p$-Muckenhoupt condition 
\begin{equation}\label{A2}
	\left (\int_Q \omega \, dx \right )
	\left (\int_{Q} \sigma \, dx\right )^{p-1} \asymp \vert Q\vert^p
\end{equation}
over all the Euclidean balls $Q\subset \R^n, \, n\geq 1 $, where  $\sigma (x)=%
\omega (x)^{1-p^{\prime}},$ it belongs to the class $A_{1}$-Muckenhoupt class if 
$$
\left ( \inf\limits_{Q} \, \omega(x) \right )^{-1} \left (\int_Q \omega \, dx \right ) \asymp \vert Q\vert.
$$
A positive measureable function $ v: \R^n\to (0,\infty)$ belongs to the  $A_{\infty}$-Muckenhoupt class if there exist  $C,\delta>0$ such that, for any ball  $Q\subset \R^n$ and its measurable  subset $E\subset Q$, the following inequality holds:
$$
\frac{v(E)}{v(Q)}\leq C\left(\frac{\vert E\vert}{\vert Q\vert}\right)^{\delta}.
$$

For  given ball $Q\subset \R^n$ and a real number $\lambda>0,$ the expression $\lambda Q,$ stands for the concentric sphere with center as $Q$ and radius $\lambda r(Q).$

A positive function $\omega: \R^n\to (0,\infty)$ satisfying the condition  
$$
\omega (2Q)\leq C\omega(Q).
$$ is called a doubling function.

For our future purposes we need to define the following Sobolev space. Let $\Omega\subset \R^N$ be a domain. Denote   $L^p(\Omega)$ the Lebesgue space of integrable functions of finite norm $\Vert f \Vert_{L^p(\Omega)} =\left ( \int\limits_{\Omega} \vert f \vert^p dz \right )^{1/p}. $ Define the Sobolev space $\mathring { \mathcal W}_{p}^1\left ( \Omega ,\omega (x)dz \right ) $ of functions with finite norm $$\Vert f \Vert _{\mathring { \mathcal W}_{p}^1(\Omega)}= \Vert f \Vert _{L^p(\Omega)}+ \Big \Vert \vert \nabla_\omega f \vert \Big  \Vert _{L^p(\Omega)} $$ where $ \vert \nabla_\omega u\vert=  \left (\omega(x)^{2/p} \vert  \nabla_{x} u \vert ^2+\vert  \nabla_{y}u\vert^2 \right )^{1/2}. $ The $\partial  / \partial x_i, \, i=1, ...n $ and $\partial / \partial y_j, \, j=1, ...m$ are partial derivatives in the distributional sense.

\smallskip

Let $p\leq q$. We will say that the pair of weights $ v\in A_\infty, \omega\in A_p ,$ defined in $R^{n},$ satisfy the condition of Chanillo-Wheeden type if
\begin{equation}\label{BC}\left ( \frac{r}{R} \right )^{1-\frac{m(n+p)}{p}\left (\frac{1}{p}-\frac{1}{q}  \right )} \left ( \frac{v \left (Q_r^x \right )}{v \left (Q_R^x \right )}\right )^{\frac{1}{q}}\leq C \left ( \frac{\omega \left (Q_r^x \right )}{\omega \left (Q_R^x \right )}\right )^{\frac{1}{p}-\frac{m}{p}\left (\frac{1}{p}-\frac{1}{q}  \right )}
\end{equation}
over all balls system $\left \{Q=Q(x,r)\right \}$ with $x\in Q_{R}^{x_0}$ and $0<r<R<\infty$
the embedding $\mathring { \mathcal W}_{p}^1 \left ( \Omega ,\omega (x)dz \right ) \rightarrow L_{q,v}(\Omega) $ is continuous for $p \leq q\leq \frac{pN}{N-p}$ and for all $ u \in  \mathring { \mathcal W}_{p}^1 \left ( \Omega ,\omega (x)dz \right ) $ 
whenever $\Omega \subset B_R^{z_0} $ in the balls of a special homogeneous space $(\R^{n+m}, \rho, dz)$ (see, \eqref{33} below).
Further, the balance condition \eqref{BC} yields compact imbedding if $1<p <q<\frac{pN}{N-p}$ and the condition \eqref{CC} below is fulfilled.

Define also the non-uniformly Sobolev space $\mathcal{\hat W} _p^1 \left ( \Omega ,\omega (x)dz \right )$ as a closure of $\mbox{Lip}_0(\Omega)$ under the norm $\Big \Vert \vert \nabla_\omega u \vert  \Big \Vert_{L^p(\Omega)} .$ Using the conditions given above concerning the function $\omega(x)$ one can show that the both spaces coincide. 

Set $E=\mathcal {\mathring W}_{p}^1 \left ( \Omega ,\omega (x)dz \right ) $ and define for all $ u \in E $ the functionals:
$$
I_{1}(u)=\frac{1}{p} \int\limits_{\Omega} \vert \nabla_\omega u\vert^p \, dz , \quad I_{2}(u)=\frac{1}{q}\int\limits_{\Omega} v(x) u(x)_+^q \, dz , \quad  I_{3}(u)=\dfrac{\mu}{\gamma}\int\limits_\Omega  u _+^{\gamma}  \, dz$$ 
$$
I(u)=I_{1} (u)-I_{2}(u)-I_{3}(u),
$$
where $ u_+(x)=\max\left \{ u(x), 0 \right \}. $
Then $I_{1}, I_{2}, I_{3}, I \in  C^1\Big ( E, \, R \Big )$ with
$$
\langle I^\prime(u), \, \phi \rangle =\langle I_{1}^\prime (u), \phi \rangle - \int\limits_{\Omega} v u_+^{q-1} \phi \, dz- \mu\int\limits_\Omega  u _+^{\gamma-1} \phi \, dz
$$
for all $u, \phi \in E. $ By definition above, $u$ is a (weak) solution of problem \eqref{1c} if $I^\prime (u)=0 $ in $ E, $ that is, for all $\phi \in E$  
\begin{equation*}\label{2}
	\int\limits_{\Omega}\vert \nabla_\omega u\vert^{p-2}\nabla_\omega u\nabla_\omega \phi \, dz 
	- \int\limits_\Omega v(x) u _+^{q-1} \phi \, dz- \mu\int\limits_\Omega  u _+^{\gamma-1} \phi \, dz=0.
\end{equation*}

In this paper we use the following Sobolev-Poincare inequalities of non-uniformly degenerated gradient.

\begin{lem} [bounded imbedding]\label{L2} Let $ q\in [p,  p(n+m)/(n+m-p],$ and let $ B_R^{z_0} $ be a fixed metric ball of the metric $\rho $ such that $\Omega \subset B_R^{z_0} $ and $m, n \geq 1.$  Assume that $v(x), \, \omega(x) $ be the functions from $A_\infty $ and $A_p $ classes of $\R^n$, respectively. Assume the condition \eqref{BC} is satisfied all over the balls system $\{ Q=Q_r^x: x \in Q_R^{x_0}, \, r<R \} . $
	
	Then it holds the inequality
	\begin{equation}\label{6a}
		\Vert f- \bar f_{v, B_R^{z_0}} \Vert_{L_{q, v}(B_R^{z_0})}\leq C_0 
		\frac{  R^{1-\frac{m(n+p)}{p}\left (\frac{1}{p}-\frac{1}{q}  \right )} v \left (Q_R^{x_0} \right )^{\frac{1}{q}}}{ \omega \left (Q_R^{x_0}  \right )^{\frac{1}{p}-\frac{m}{p}\left (\frac{1}{p}-\frac{1}{q}  \right )}}\Vert \nabla_\omega f \Vert_{L_p(B_R^{z_0})} 
	\end{equation}
	with $\vert \nabla_\omega f \vert^p=\left (\omega(x)^{2/p}\vert \nabla_x f \vert ^2+\vert \nabla_y  f\vert ^2 \right  )^{p/2} $
	for any Lipshitsz continues function $u$ in the $\rho$ metric ball $ B_R^{x_0}, $ where $\bar f_{v, B_R^{z_0}}=\frac{1}{v(B_R^{z_0})} \int\limits_{B_R^{z_0}} v f \, dx , $ the constant $C_0$ depends on the constants $C,\delta$ from the $A_p $ condition and $n, q.$
\end{lem}

For a bounded convex domain $\Omega$ of weight $\omega \in A_{p}$, using Lemma 1, we get 
$$ \mathring  { \mathcal W}_{p}^1 \left ( \Omega ,\omega (x)dz \right ) \subset  \mathring  { \mathcal W}_{1}^1 (\Omega)\subset \subset L_\gamma(\Omega )   \quad 1<\gamma <N/(N-1).$$

\smallskip

We use the next result from \cite{FM5} to prove Lemma \ref{L2}.  Let the functions $v, \omega_1, \omega_2, ...\omega_n, : \, \R^n \to (0, \infty)$ are measurable functions. 
\begin{teo}. \label{tPS}
	Let \,  $ q\geq p \ge 1, \quad  \left( \R^N,  \rho, dz \right) $ be a homogeneous space and assume that all balls are convex. Let $ B_R^{z_0} \subset \R^N $ be a fixed metric ball. Suppose that $ v\in A_\infty,  \, \,  \omega_i \in A_p$, for $i=1,2,...N$, and that
	\begin{equation} \label {BAC}
		l_i(B^*)\vert B \vert^{-(\frac{1}{p}-\frac{1}{q})} \Big( \frac{1}{\vert  B \vert}\int\limits_{B\cap \Omega} v \, dz \Big) ^{\frac{1}{q}} \Big( \frac{1}{\vert  B \vert}\int\limits_{B\cap \Omega} \omega_i^{1-p^\prime} \, dz \Big)^{\frac{1}{p^\prime}}\leq A, \quad i=1,2,...N,
	\end{equation}
	all over the $\rho$-quasimetric balls $B=B(z, r)$ having center in \, $ z\in B_R^{z_0}$ and radius \, $r<R $, where $l_i(B^*)=\sup \, \Big \{\vert x_i-y_i  \vert : \, x,y \in B^* \Big \} $ and $B^*=B\left (z, 5K_0^2 \, r \right )$. Then, for all Lipschitz continuous functions $f:B_R^{z_0} \to \R $  the Poincare's type inequality
	\begin{equation} \label {PS}
		\Big ( \int \limits_{B_R^{z_0}} \vert f -\bar f_{v, B_R^{z_0}}\vert ^q v \, dz \Big)^{\frac{1}{q}}\leq  {C_0} A \,  \sum \limits_{i=1}^N
		\Big ( \int\limits_{B_R^{z_0}} \vert  f_{z_i} \vert ^p \omega_i \, dz \Big)^{\frac{1}{p}}, \quad f\in C_0^1(\Omega).
	\end{equation}
	holds with a constant $C_0$ depending only on $p,q, n$ and on the Muckenhoupt condition constants $C, \delta$.
\end{teo}

For
$ \omega :\,\R^n \to  (0,\infty )$ denote $\sigma (x)=%
\omega (x)^{1-p^{\prime}}$  is in the Muckenhoupt's $A_{p}$-class all over the $n$-dimensional  Euclidean balls.
Define a function $h_{x}(.):[0,\infty )\to (0,R)$ as
\[
h_x(t)=t\left(\frac{1}{\vert Q(x,t) \vert}\int\limits_{Q(x,t)}\, \sigma(s)\, ds \right)^{\frac{1}{p^{,}}},\quad x\in \R^n,
\]
and assume that  $h_{x}(0)=0$,  $\displaystyle\lim_{t\to + \infty} h_x(t)=+\infty$ for a fixed $x\in \R^n $.  Then we may consider an inverse function $h_{x}^{-1}(.):[0, R)
\to  (0,\infty )$ defined as
\begin{equation*}
	h_{x}^{-1}(v)=\inf \Big \{\rho >0:\,\,h_{x}(\rho )\geq v\Big \},\quad v>0.
\end{equation*}

We can define a quasi-metric $\rho$ on $\mathbb{R}
^{N}=
\mathbb{R}
^{n}\times
\mathbb{R}
^{m}=\{z=(x,y) \>| x \in \mathbb{R}
^{n}\>, y\in
\mathbb{R}
^{m} \}$ as follows: for any $z_{1}=(x_{1},y_{1}), z_{2}=(x_{2},y_{2}) \in \mathbb{R}
^{N}$ we put
\begin{equation}
	\rho (z_{1},z_{2})=\max \Big \{|x_{1}-x_{2}|,\,%
	\,h_{x_{1}}^{-1}(|y_{2}-y_{1}|),\,\,h_{x_{2}}^{-1}(|y_{2}-y_{1}|)\Big \}.
	\label{33}
\end{equation}
It is not difficult to check that $\rho :
\mathbb{R}^{N}\times
\mathbb{R}
^{N}\rightarrow \lbrack 0,\infty )$ is a quasi-metric satisfying the triangle inequality \begin{equation}
	\rho (z_{1},z_{2})\leq K_0\Big (\rho (z_{1},z_{3})+\rho (z_{2},z_{3})\Big)
	\label{35}
\end{equation}%
with a constant $K_0\geq 1$ independent of $z_{1},z_{2},z_{3}\in
\mathbb{R}
^{N},$ (see, e.g., \cite{HA}). 
Therefore, the constructed quasi-metric space endowed with the Lebesgue measure is a homogeneous space. Set $\left (\R^{n+m}, dx, \rho\right )$ for it.

Now we give the proof of inequality \eqref{6a} by using Theorem \ref{tPS}  (see, \cite{FM5}).We also need the cited above homogeneous space with metric \eqref {33} and measure $dx.$ Let us set $ z=(x, y): \, x\in \R^n, \, y\in \R^m,  \, v(x, y):=v(x), \omega_1(x, y)=\omega_2(x, y)=, ...,=\omega_n(x, y):=\omega(x)$ and $\omega_{n+1}(x, y)=\omega_{n+2}(x, y)=, ..., =\omega_{n+m}(x, y):=1.$
From inequality \eqref{PS} it follows that 
$$
\Vert u- \bar u_{v, B_R^{z_0}} \Vert_{L_{q, v}(B_R^{z_0})}\leq C_0 A \Big (\iint\limits_{B_R^{z_0}}\bigg ( \omega^{\frac{2}{p}} \vert  \nabla_xu \vert^2+\vert \nabla_y u \vert ^2\bigg )^\frac{p}{2} \, dx dy \Big )^{\frac{1}{p}}.
$$
In condition \eqref{BAC}, we can take the supremum of 
\begin{equation}\label{At}
	C r \left [ r^{n+m} \left (\frac{1}{r^n} \int_{Q_r^x} \sigma \, dx \right )^{\frac{m}{p^{,}}}       \right ]^{-\left (\frac{1}{p}-\frac{1}{q} \right )}
	\left ( \frac{1}{r^n} \int_{Q_r^x} v \, dx \right )^{\frac{1}{q}} \, \left ( \frac{1}{r^n} \int_{Q_r^x} \sigma \, dx \right )^{\frac{1}{p^{,}}}
\end{equation}
all over the balls system $\left \{Q=Q_r^x: \, x\in Q_R^{x_0}, \,  r<R \right \}.$

The expression \eqref{At} becomes

$$
C r^{1-(n+m)\left(\frac{1}{p}-\frac{1}{q}  \right)+n\left(\frac{1}{p}-\frac{1}{q}  \right)-\frac{nm}{p}\left(\frac{1}{p}-\frac{1}{q}\right )} \left(  \int_{Q_r^x} v \, dx \right)^{\frac{1}{q}} \left ( \int_{Q_r^x} \sigma \, dx \right )^{\frac{1}{p^{,}}\left(1-m\left(\frac{1}{p}-\frac{1}{q}  \right )\right)\frac{p-1}{p-1}}
$$
Using the condition $A_p$
$$
\left ( \frac{1}{r^n} \int_{Q_r^x} \sigma(s) \, ds \right )^{(p-1)}\leq C  \left ( \frac{1}{r^n} \int_{Q_r^x} \omega (s) \, ds \right )^{-1}
$$ and assuming $ \frac{1}{p}-\frac{m}{p}\left (\frac{1}{p}-\frac{1}{q}  \right ) \geq 0, $
\eqref{At} can be estimated by the expression
$$
C r^{1-(n+m)\left (\frac{1}{p}-\frac{1}{q}  \right  )+n\left(\frac{1}{p}-\frac{1}{q}  \right )-\frac{nm}{p}\left(\frac{1}{p}-\frac{1}{q}\right )} \left (  \int_{Q_r^x} v \, dx \right )^{\frac{1}{q}} \left ( \int_{Q_r^x} \omega \, dx \right )^{\frac{1}{p^{,}}\left(-1+m\left(\frac{1}{p}-\frac{1}{q}  \right )\right)\frac{1}{p-1}}
$$
$$
=C\,  r^{1-(n+m)\left (\frac{1}{p}-\frac{1}{q}  \right )} \, r^{n\left (\frac{1}{p}-\frac{1}{q} \right )}  \, r^{-\frac{nm}{2}\left (\frac{1}{p}-\frac{1}{q} \right ) \left (2-p\right )}  \left ( \int_{Q_r^x} v \, dx \right )^{\frac{1}{q}} \left ( \int_{Q_r^x} \omega \, dx \right )^{-\frac{1}{p}+\frac{m}{2}\left (\frac{1}{p}-\frac{1}{q}  \right )}
$$
\begin{equation}\label{l21}
	=\frac{  r^{1-\frac{m(n+p)}{p}\left (\frac{1}{p}-\frac{1}{q}  \right )} v \left (Q_r^{x_0} \right )^{\frac{1}{q}}}{ \omega \left (Q_r^{x_0}  \right )^{\frac{1}{p}-\frac{m}{p}\left (\frac{1}{p}-\frac{1}{q}  \right )}}
\end{equation}
Now, it remains to demand the condition \eqref{BC} in order to get an upper bound for the constant $A$; in conclusion, the right-hand side of \eqref{l21} will be less than or equal to
\begin{equation}\label{bcon}
	\frac{  R^{1-\frac{m(n+p)}{p}\left (\frac{1}{p}-\frac{1}{q}  \right )} v \left (Q_R^{x_0} \right )^{\frac{1}{q}}}{ \omega \left (Q_R^{x_0}  \right )^{\frac{1}{p}-\frac{m}{p}\left (\frac{1}{p}-\frac{1}{q}  \right )}}
\end{equation}
(i.e., we can set the constant $A$ equal to this ratio).

Note that \eqref{bcon}  is exactly equal to
$$
C_{0}R\frac{v(B^{z_{0}}_{R})^{\frac{1}{q}}}{\omega(B^{z_{0}}_{R})^{\frac{1}{p}}}.
$$

Therefore, the following inequality holds:
$$
\Big(\frac{1}{v(B^{z_{0}}_{R})}\int_{B^{z_{0}}_{R}}^{ }\vert f- \bar f_{v, B_R^{z_0}} \vert ^{q}vdz\Big)^{1/q}\leq C_{0}R\Big(\frac{1}{\omega(B^{z_{0}}_{R})}\int_{B^{z_{0}}_{R}}^{ }\vert \nabla_{\omega}f \vert ^{p}dz\Big)^{1/p}.
$$


\medskip

\begin{lem}[compact imbedding]\label{L3} Let $ q\in (0,  p(n+m)/(n+m-p),$and let $ B_R^{z_0}\subset \R^{n+m} $ be a fixed metric ball of metric $\rho $ from \eqref{33} such that $\Omega \subset B_R^{z_0} $ and $m, n \geq 1.$   Assume that $v(x), \omega(x): Q_{2R}^{x_0} \to (0,\infty)$ be the functions from $A_\infty $ and $A_p$ classes of $\R^n$, respectively,  and satisfying \eqref{BC}
	all over the balls system $\left \{Q=Q(x,r): \, \,   \text{with} \, \,  x\in Q_R^{x_0}\, \,  \text{ and} \, \,  r<R \right \}. $

	If 
	\begin{equation}\label{CC}
		\frac{ {C} R^{1-\frac{m(n+p)}{p}\left (\frac{1}{p}-\frac{1}{q}  \right )} v \left (Q_R^{x_0} \right )^{\frac{1}{q}}}{ \omega \left (Q_R^{x_0}  \right )^{\frac{1}{p}-\frac{m}{p}\left (\frac{1}{p}-\frac{1}{q}  \right )}}\to 0 \quad \text{as}\quad r\to 0 \quad \text{uniformly on} \quad x \in Q_R^{x_0} 
	\end{equation}
	then the set of Lipshitsz continuos functions $u(z): B_R^{z_0}\subset \R^{n+m} \to \R $ with 
	$$
	\int\limits_{\Omega}\bigg ( \omega^{2/p} \vert  \nabla_xu \vert^2+\vert \nabla_y u \vert ^2\bigg )^{p/2} \, dz  \leq M 
	$$
	is compactly imbedded into $L_{q,v}(\Omega).$
\end{lem}

The assertion can be proved by using Lemma \ref{L2} as it was done in \cite[Lemma 4]{FC}.

\bigskip

\section{Main results}

\begin{teo}.  \label{t1}
	Let $\Omega \subset \R^{n+m} $ be a bounded domain, $ p<q<p(n+m)/(n+m-p)$ and $1<\gamma<N/(N-1).$ Assume the conditions $\omega \in A_p, \, v\in A_\infty,  $ \, $\mu \in (0, \Lambda),$ $\Lambda$ is sufficiently small, and  \eqref{BC}, \eqref{CC} are fulfilled.  Then the problem \eqref{1c} has two different positive solutions from space $\mathring  { \mathcal W}_{p}^1 \left (\Omega ,\omega (x)dz \right )$. 
\end{teo}

\bigskip

\begin{rem}. The range condition $1<\gamma<N/(N-1)$ in Theorem \ref {t1} is assumed to have the inclusion $ \mathring  { \mathcal W}_{p}^1 \left ( \Omega ,\omega (x)dz \right ) \subset  \subset L_\gamma (\Omega)$. In this theorem, we may consider the range $1<\gamma<p$ if we add the condition $ v^{-\gamma/(q-\gamma)}\in L_{1, loc}$ to the list of assumptions,  which provide the inclusions $ \mathring  { \mathcal W}_{p}^1 \left ( \Omega ,\omega (x)dz \right ) \subset \subset L_{q, v}(\Omega) \subset L_\gamma(\Omega ).$   
\end{rem}

Applying Theorem 4 to the case $v=a\omega (x), a>0,$ we obtain the following corollary (see e.g.\cite{FC}).

\begin{cor}.
	Let $\Omega \subset R^{N}$ be a bounded domain, let $1<\gamma <N/(N-1),$ let the function $\omega$ of the variable x belongs to the class $A_{p}$ in $R^{n}$, let $\mu \in (0,\Lambda),$ where $\Lambda$ is small enough, and let the following assumptions holds:
	$$
	q \in \left(p, \frac{pm(n+p)}{m(n+p)-p^{2}}\right) \quad for \quad m>p,
	$$
	$$
	q \in \left(p, \frac{pn(p^{2}-p+1)+p^{2}}{n(p^{2}-p+1)+p-p^{2}}\right) \quad for \quad m=1,
	$$
	$$
	q \in \left(p, \frac{p(n+p)}{n}\right) \quad for \quad m=p.
	$$
	Then the following problem:
	$$
	\nabla_\omega( \vert \nabla_\omega u\vert^{p-2}\nabla_\omega u)+a\omega (x) u^{q-1}+ \mu u^{\gamma-1}=0\quad \text{in} \quad \Omega, 
	$$
	$$
	u>0 \quad \text{in} \quad \Omega, \\
	$$
	$$
	u=0 \quad \text{in} \quad \partial \Omega
	$$	
	for $\alpha >0$ has at least two different positive solutions from the space $\mathring { \mathcal W}_{p}^1 \left ( \Omega ,\omega (x)dz \right ).$
\end{cor}

\bigskip

\textit{Essential properties}.
Set $ E=\mathring  { \mathcal W}_{p}^1 \left ( \Omega ,\omega (x)dz \right )$ be a Banach space of functions which will be used essentially. It is possible to solve problem \eqref{1c}, by finding critical points of the following functional.
$$
I(u)=\frac{1}{p}\int\limits_{\Omega} \vert \nabla_\omega u\vert^p  \, dz -\frac{1}{q} \int \limits_{\Omega}v u_+^q dz - \dfrac{\mu}{\gamma}\int\limits_\Omega  u _+^{\gamma}  \, dz.
$$
This functional has continuous derivative in the sense of Gateaux: $I^{'}(u)\in E^{*} $ from the dual space, the derivative maps $u\in E\mapsto I^{'}(u)\in E^{*} $ as follows:  $\forall \phi \in \mathring  { \mathcal W}_{p}^1 \left ( \Omega ,\omega (x)dz \right )$
$$
\left \langle I^\prime(u), \phi \right \rangle =\int\limits_{\Omega} \vert \nabla_\omega u\vert^{p-2}\nabla_\omega u\nabla_\omega \phi  \, dz 
- \int\limits_\Omega v(x) u _+^{q-1} \phi \, dz-\mu\int\limits_\Omega  u _+^{\gamma-1} \phi \, dz.
$$
Denote the first functional 
$$
I_{1}(u)=\frac{1}{p}\int\limits_{\Omega} \vert \nabla_\omega u\vert^p  \, dz 
$$
and confirm that it belongs to $C^1(E,R). $ 

\textit{The functional $I_{1}(u)\in C^1(E, \R). $ } Take $\forall \phi \in \mathring  { \mathcal W}_{p}^1 \left ( \Omega ,\omega (x)dz \right )$
$$\left \langle I_{1} ^\prime(u), \phi \right \rangle =\int\limits_{\Omega} \vert \nabla_\omega u\vert^{p-2}\nabla_\omega u\nabla_\omega \phi  \, dz .$$
Let any sequence $u_k \to u$  in $\mathring  { \mathcal W}_{p}^1 \left ( \Omega ,\omega (x)dz \right ),$ getting 
$$\vert \left \langle  I_{1} ^\prime(u_k)-I_{1} ^\prime(u), \phi \right \rangle \vert=\vert\int\limits_{\Omega}\left[\vert \nabla_\omega u_k\vert^{p-2}\nabla_\omega u_k-\vert \nabla_\omega u\vert^{p-2}\nabla_\omega u \right]\nabla_\omega \phi \, dz\vert $$
in case of $ 1\leq p\leq 2 $ (see \cite{FAMM}, Lemma 1.6)
$$\leq C(p)\int\limits_{\Omega}\vert\nabla_\omega u_k-\nabla_\omega u\vert^{p-1}\nabla_\omega \phi  \, dz \leq  C(p)\Vert u_k-u\Vert_{\mathring  { \mathcal W}_{p}^1 \left ( \Omega ,\omega (x)dz \right )}^{p-1}\Vert\phi \Vert_{\mathring  { \mathcal W}_{p}^1 \left ( \Omega ,\omega (x)dz \right )},$$
in case of $ p>2 $
$$\leq(p-1)\int\limits_{\Omega}\vert\nabla_\omega u_k-\nabla_\omega u\vert \left( \vert \nabla_\omega u_k\vert^{p-2}+ \vert \nabla_\omega u\vert^{p-2} \right)\vert\nabla_\omega \phi \vert \, dz$$
$$=(p-1)\int\limits_{\Omega}\vert\nabla_\omega u_k-\nabla_\omega u\vert \left( \vert \nabla_\omega( u_k-u)+\nabla_\omega u\vert^{p-2}+ \vert \nabla_\omega u\vert^{p-2} \right)\vert\nabla_\omega \phi \vert  \, dz$$
$$\leq(p-1)\int\limits_{\Omega}\vert\nabla_\omega u_k-\nabla_\omega u\vert \left( \vert \nabla_\omega (u_k-u)\vert^{p-2}+ 2\vert \nabla_\omega u\vert^{p-2} \right)\vert \nabla_\omega \phi \vert \, dz$$
$$=(p-1)\int\limits_{\Omega}\vert\nabla_\omega u_k-\nabla_\omega u\vert^{p-1} + 2\vert \nabla_\omega u\vert^{p-2}\vert\nabla_\omega u_k-\nabla_\omega u\vert \vert \nabla_\omega \phi \vert \, dz$$
$$\leq(p-1)\Vert u_k-u\Vert_{\mathring   { \mathcal W}_{p}^1 \left ( \Omega ,\omega (x)dz \right )} ^{p-1} \Vert\phi \Vert_{\mathring   { \mathcal W}_{p}^1 \left ( \Omega ,\omega (x)dz \right )}+ 2(p-1)\int\limits_{\Omega}\vert \nabla_\omega u\vert^{p-2}\vert\nabla_\omega u_k-\nabla_\omega u\vert \vert \nabla_\omega \phi \vert \, dz$$
$$\leq(p-1)\Vert u_k-u\Vert_{\mathring   { \mathcal W}_{p}^1 \left ( \Omega ,\omega (x)dz \right )} ^{p-1} \Vert\phi \Vert_{ \mathring{ \mathcal W}_{p}^1 \left ( \Omega ,\omega (x)dz \right )}$$
$$+ 2(p-1)\Vert u_k-u\Vert_{\mathring   { \mathcal W}_{p}^1 \left ( \Omega ,\omega (x)dz \right )}^{p-1}  \left(\int\limits_{\Omega}\vert\nabla_\omega u\vert^{\frac{p(p-2)}{p-1}}\vert \nabla_\omega \phi \vert^{\frac {p}{p-1}} \, dz \right)^\frac{p-1}{p} ,$$

$$\int\limits_{\Omega}\vert\nabla_\omega u\vert^{\frac{p(p-2)}{p-1}}\vert \nabla_\omega \phi \vert^{\frac {p}{p-1}} \, dz \leq \left(\int\limits_{\Omega}\vert\nabla_\omega u\vert^{\frac{p(p-2)}{p-1}\frac{p-1}{p-2}} \, dz \right)^{\frac{p-2}{p-1}}\left(\int\limits_{\Omega}\vert\nabla_\omega \phi \vert^{\frac{p}{p-1}(p-1)} \, dz\right)^{\frac{1}{p-1}} $$
$$\leq\Vert u\Vert_{\mathring { \mathcal W}_{p}^1 \left ( \Omega ,\omega (x)dz \right )}^{\frac{p(p-2)}{p-1}}\Vert\phi \Vert_{\mathring { \mathcal W}_{p}^1 \left ( \Omega ,\omega (x)dz \right )}^{\frac{p}{p-1}}.$$
So it shows tbe continuity of the functional.

\hspace{10cm} $\Box $

\bigskip

\textit{ The functional  $ I_{2}(u)=\frac{1}{q}\int\limits_{\Omega} v(x) u(x)_+^q \, dx  \in C^1(E, E^*) .$}  Then 
$$
\left \langle  I_{2} ^\prime (u), v\right \rangle =\int\limits _{\Omega}v(x) u_+(x)^{q-1} v(x) \, dx , \quad \forall u, v \in \mathring  { \mathcal W}_{p}^1 \left ( \Omega ,\omega (x)dz \right ).
$$

Take $q>1$ set the function $F(t)=t^q$, where $t\geq 0$ and $F(t)=0$, where $t<0.$ Clearly 
$F(t)\in C^1(\R)$ and $F^\prime(t)=qt^{q-1}$, $t\geq 0$ and  $F^\prime(t)=0$, $t<0, $ and  
$ u_+(x)=F\left (u(x)\right ). $ Using the above notations for $ \forall u, \phi \in \mathring  { \mathcal W}_{p}^1 \left ( \Omega ,\omega (x)dz \right )$ and $t\in \R\setminus\{ 0 \}$ we have
$$
\frac{1}{t} \left [ I_{2} (u+t\phi)-I_{2} (u) \right ]=\frac{1}{t}\int\limits_{\Omega} v(x)\left [ F(u(x)+t\phi(x))-F(u(x))\right ] \, dx 
$$
$$
=\int\limits_{\Omega} v(x)F^\prime \left (u(x)+\theta \phi(x)\right ) \, dx =\int\limits_{\{ x: u+\theta \phi>0\}} v \left ( u+\theta \phi \right )^{q-1} \phi(x) dx
$$
$$
=\int\limits_{\Omega} v(x) \left ( u(x)+\theta \phi(x) \right )^{q-1}\phi (x) \,   \chi_{\{u+\theta \phi>0\}} \, dx  
$$
$$
\to \int\limits_{\Omega} v(x)  u(x)^{q-1}\phi (x) \,   \chi_{\{u>0 \}} \, dx\quad \text{as}\quad t\to 0
$$
on base of Lebesgue majorant theorem since 
$$
v(x) \left ( u(x)+\theta \phi(x) \right )^{q-1}\phi (x) \,   \chi_{\{u+\theta \phi>0\}} \to v u_+^{q-1}\phi \quad\text{as}\quad t\to 0 $$ 
a.e. $ x\in \Omega.$
There exists an integrable majorant function, by the condition \eqref{BC} and imbedding $\mathring  { \mathcal W}_{p}^1 \left ( \Omega ,\omega (x)dz \right ) \subset   L_{q, v} (\Omega)  : $
$$
v(x) \left ( u(x)+\theta \phi(x) \right )^{q-1} \vert \phi (x)\vert  \,   \chi_{\{u+\theta \phi>0\}}\leq  v \left (\vert u\vert^{q-1} \vert \phi \vert+\vert \phi \vert^q\right ) \in L_1(\Omega),
$$
by using Holder's inequality,
$$
\int\limits_{\Omega} v \vert u\vert^{q-1} \vert \phi \vert \leq \Vert u \Vert_{L_{q, v}(\Omega)}^{q-1}\Vert \phi \Vert_{L_{q, v}(\Omega)}\leq C \Vert u \Vert _{\mathring { \mathcal W}_{p}^1 \left ( \Omega ,\omega (x)dz \right )}^{q-1}\, \Vert \phi \Vert _{\mathring  { \mathcal W}_{p}^1 \left ( \Omega ,\omega (x)dz \right )}^{q-1}<\infty  \cdot
$$

\hspace{10cm} $\Box$

\textit {The functional $ I(u)\in C^1(E, \R) .$} 

On basis of the property $I^\prime(u): E\to E^*$ and $\forall u, \phi \in \mathring  { \mathcal W}_{p}^1 \left ( \Omega ,\omega (x)dz \right )$ we have 
$$
\left \langle I^\prime(u), \phi \right \rangle =\int\limits_{\Omega}\vert \nabla_\omega u \vert^{p-2}\nabla_\omega u\ \nabla_\omega \phi  \, dz- \int\limits_\Omega v(x) u_+^{q-1} \phi(z )\, dz-\mu\int\limits_\Omega  u _+^{\gamma-1} \phi(z) \, dz. 
$$

\hspace{10cm} $\Box $

The functional $I(u)$ is lower semi-continuous (see e.g. \cite{FC}).

\smallskip

\proof of Theorem \ref{t1}.

\textit{Palais-Smale condition}. In order to show PS -condition we should demonstrate the sequence $ \{u_n\}\in E $ is compact (set $ E=\mathring { \mathcal W}_{p}^1 \left ( \Omega ,\omega (x)dz \right ) $), i.e. contains a convergent subsequence $ u_{n_k} \to u\in E. $

Firstly, establish the boundedness of $\{u_{n}\} $
in $ E $. Using 1) it follows
$$
I_{1}(u_n)- I_{2}(u_n)- I_{3}(u_n) \leq M, 
$$
then
\begin{equation}\label{10a}
	I_{1}(u_n)\leq  I_{2}(u_n)+I_{3}(u_n) + M.
\end{equation}

On other hand, using condition 2), $\Vert I^{\prime }(u_{n})\Vert_{E^*}=o(1) $ as $ n\rightarrow
\infty $ it follows
\begin{equation}\label{1.6}
	\left \langle I^{\prime }(u_n),\phi \right \rangle
	=o(1)\, \left\Vert \phi \right\Vert_E , \, \, \forall v \in E .
\end{equation}
In particularly, setting $ \phi=u_n, $ we get
\[
\left\langle I^{\prime }(u_n),u_n\right\rangle
=o(1)\, \left\Vert u_{n}\right\Vert_E
\] i.e.

$$
\left \langle I^\prime(u_n), u_n \right \rangle =\int\limits_{\Omega}\vert \nabla_\omega u_n \vert^{p-2}\nabla_\omega u_n\ \nabla_\omega u_n  \, dz- \int\limits_\Omega v(x) (u_n)_+^{q-1} u_n \, dz$$$$-\mu\int\limits_\Omega  (u_n)_{+} ^{\gamma-1} u_n \, dz 
=o(1)\, \left\Vert u_{n}\right\Vert_E.$$

By utilizing the identities
$$
v(u_n)_{+}^{q-1}u_n \equiv v(u_n)_{+}^{q},\quad \mu (u_n)_{+}^{\gamma-1}u_n\equiv \mu (u_n)_{+}^{\gamma},
$$
we get
\begin{equation}\label{YC}
	\frac{1}{q}\int\limits_\Omega v(x) (u_n)_+^{q-1} u_n(z )\, dz+\frac{\mu}{q}\int\limits_\Omega  (u_n)^{\gamma-1} u_n \, dz $$$$=o(1)\, \left\Vert u_{n}\right\Vert_E+\frac{1}{q}\int\limits_{\Omega}\vert \nabla_\omega u_n \vert^{p-2}\nabla_\omega u_n\ \nabla_\omega u_n  \, dz$$ $$
\end{equation}
From \eqref{YC}, \eqref {10a} and $\mu\geq 0$, we conclude
\begin{equation*}
	\left( \frac{1}{p}-\frac{1}{q}\right) \Lambda(u_n)
	\leq M+o(1)\left\Vert u_{n}\right\Vert_E .
\end{equation*}
From this for $ q>p, $ it follows
\begin{equation*}
	\Vert u_n\Vert_{E}^p\leq \frac{pqM}{
		q-p}+o(1)\left \Vert u_{n}\right\Vert_E.
\end{equation*}
Connecting this with Young's inequality, it follows
\begin{equation}\label{2.7}
	\Vert u_{n}\Vert _E\leq C(M).
\end{equation}

This completes the boundedness of the sequence $ \{u_{n}\} $ in $ E $.

Applying well-known
fact, there exists a weak convergent subsequence $ u_{n_{k}} \rightarrow u $ in
$ E $. Denote it again $ u_{n}. $ Utilizing the compact embedding Lemma \ref{L3}, we get the strong convergence $ u_{n}\rightarrow u\ $ in $  L_{q, v}(\Omega) $ , i.e.%
\begin{equation*}
	\Vert u_{n}-u\Vert _{L_{q, v}(\Omega)} \rightarrow 0
\end{equation*}
Now, we are ready to show the strong convergence $ u_{n}\rightarrow u $ in $ E $%
. For this, taking $ \phi =u_{n}-u $ in \eqref{1.6}:

$$\int\limits_{\Omega}\vert \nabla_\omega u_n \vert^{p-2}\nabla_\omega u_n\ \nabla_\omega (u_n(z)-u(z)) \, dz$$

\begin{equation}\label{13a}
	\begin{split}	
		-\int\limits_\Omega v(x) (u_n)_+^{q-1}\left ( u_n(z)-u(z) \right )\, dz=o(1)\left\Vert u_{n}-u\right\Vert_E .
	\end{split}
\end{equation}%
From this, since $ u_{n}\rightarrow u $ in $ L_{q, v}(\Omega)$ and using Holder's inequality, it follows
\[
\Big \vert \int\limits_\Omega v(x) (u_n)_+^{q-1}\left ( u_n(x)-u(x) \right )\, dx \Big \vert  \leq \left\Vert u_{n}-u\right\Vert _{L_{q, v}(\Omega)} \Vert (u_n)_+\Vert_{L_{q, v}(\Omega)}^{q-1}
\] 
\[
=C(M)  \Vert u_n-u\Vert_{L_{q, v}(\Omega)}  \to 0 \quad \text{as} \quad n\to \infty.
\]
where also has been used Lemma \ref{L2} and the estimate \eqref{2.7}, in order to obtain the boundedness $ \{u_{n}\} $ in $ L _{q,v}(\Omega). $

Also from \eqref{2.7} and the assumption  $\omega^{1-p^{,}}\in L_{1, loc}(\R^n)$ using the Holder inequality, it follows 
\begin{equation}\label{Hol}
	\int_{\Omega}\left ( \vert \nabla _x u_n \vert+\vert \nabla_y u_n \vert \right ) \, dz \leq \left ( 
	\int_{\Omega} \Big( \omega^{\frac{2}{p}}\vert \nabla _x u_n \vert^2+\vert \nabla_y u_n \vert ^2 \Big )^{\frac{p}{2}} \, dz \right )^{1/p} \times 
\end{equation}$$ 
\times \left ( 
\int _{\Omega}\left (\omega^{1-p^{\prime}} +1 \right )  \, dz \right )^{1/p^{\prime}}\leq C_2(M). $$
Therefore, using the compact imbedding $\mathring{W}_1^1(\Omega)\subset \subset L_{r}(\Omega)$ for any $r\in \left [1, \,  N/(N-1)\right )$ there exists a subsequence $u_{n_k}$ (denote the subsequence like $u_n$) such that $u_n \to u$ strongly in $L_{r}(\Omega).$ On the basis of \eqref{YC}, we assume that, $\Vert u_n \Vert_{L_\gamma (\Omega)} \leq C_3(M)$ therefore, and using Holder's inequality for $1\leq \gamma <N/(N-1)$ we have 
\[
\Big \vert \int\limits_\Omega  (u_n)_+^{\gamma -1}\left ( u_n(z)-u(z) \right )\, dz \Big \vert  \leq \left\Vert u_{n}-u\right\Vert _{L_{\gamma}(\Omega)} \Vert (u_n)_+\Vert_{L_{\gamma}(\Omega)}^{\gamma -1}
\] 
\[
\leq C_3(M)  \Vert u_n-u\Vert_{L_{\gamma}(\Omega)} =o(1) \quad \text{as} \quad n\to \infty.
\]
Thus, from \eqref{13a}
\begin {equation}\label{14a}\int\limits_{\Omega}\vert \nabla_\omega u_n \vert^{p-2}\nabla_\omega u_n\ \nabla_\omega (u_n (z)-u(z)) \, dz= o(1)+o(1)\left\Vert u_{n}-u\right\Vert_E .
\end{equation}
or
$$
\int\limits_{\Omega}\vert \nabla_\omega u_n \vert^{p-2}\nabla_\omega u_n\ \nabla_\omega (u_n (z)-u(z)) \, dz+ 
\int\limits_{\Omega}\vert \nabla_\omega u \vert^{p-2}\nabla_\omega u\ \nabla_\omega (u_n (z)-u(z)) \, dz
$$
$$=\int\limits_{\Omega}\vert \nabla_\omega u \vert^{p-2}\nabla_\omega u\ \nabla_\omega (u_n (z)-u(z)) \, dz+ o(1)+o(1)\left\Vert u_{n}-u\right\Vert_E .
$$
We have

$$\int\limits_{\Omega}\vert \nabla_\omega (u_n-u) \vert^{p} \, dz\leq \int\limits_{\Omega}\vert \nabla_\omega u \vert^{p-2}\nabla_\omega u\ \nabla_\omega (u_n (z)-u(z)) \, dz+ o(1)+o(1)\left\Vert u_{n}-u\right\Vert_E$$

The first term in the right hand side tends to zero since $u_n \to u $ weakly in $\mathring { \mathcal W}_{p}^1 \left ( \Omega ,\omega (x)dz \right ):$ 
\[
\int\limits_{\Omega}\vert \nabla_\omega u \vert^{p-2}\nabla_\omega u\ \nabla_\omega (u_n (z)-u(z)) \, dz\to 0 \quad \textit{as} \quad k \to \infty.
\]
connecting this with \eqref{14a}, infers
\[
\Vert u_n-u\Vert _E^p=o(1)+o(1)\left\Vert u_{n}-u\right\Vert_E .
\]
Using here
Young's inequality, it follows
\[
\Vert u_n-u\Vert _E=o(1) \quad \text{as} \quad k\to \infty.
\]
Therefore, $ u_{n}  \to  u $ in $ E, $ which means the sequence $\{ u_n\}$ is compact in E.$ 
$
This proves PS -property for the functional $I(u).$

\hspace{10cm} $\Box$

\textit{Existence of the first solution.}
we apply the MPT in order to show the existence of solution
for the problem \eqref{eq1}.

For $ \Vert u \Vert_E \leq 1 $ we have%
$$
I(u)=\frac{1}{p}\int\limits_{\Omega} \vert \nabla_\omega u\vert^p  \, dz  -\frac{1}{q}\int\limits_\Omega v(x) u_+^q(z) \, dz- \frac{\mu}{\gamma}\int\limits_\Omega  u_{+}^{\gamma}(z) \, dz.
$$
Since domain is bounded we may assume that $\Omega \subset B_R^{z_0} $ . 
Using the Sobolev type inequality (Lemma \ref{L2})

$$
\Vert u \Vert_{L_{q, v}(Q_R^{x_0})} \leq C_0 
\frac{  R^{1-\frac{m(n+p)}{p}\left (\frac{1}{p}-\frac{1}{q}  \right )} v \left (Q_R^{x_0} \right )^{\frac{1}{q}}}{ \omega \left (Q_R^{x_0}  \right )^{\frac{1}{p}-\frac{m}{p}\left (\frac{1}{p}-\frac{1}{q}  \right )}}\Vert \nabla_\omega u \Vert_{L_{p}(Q_R^{x_0})}
$$
and

$$
\Vert u \Vert_{L_{\gamma, v}(\Omega)} \leq C\vert \Omega \vert^{\frac{1}{\gamma}-\frac{1}{N^{,}}}\Vert\nabla_z u\Vert_{L_{1}(\Omega)} \leq C\vert \Omega \vert^{\frac{1}{\gamma}-\frac{1}{N^{,}}}\left( \int\limits_\Omega (1+\omega^\frac{-p^{,}}{p})dz\right)^\frac{1}{p^{,}}\Vert u\Vert_{E} 
$$
for $ \mu \in (1,N/(N-1))$,
we get 
$$
I(u)\geq \frac{1}{2p}\Vert u \Vert_E^p +\frac{1}{4p}\Vert u \Vert_E^p -\frac{1}{q}\Bigg (C_0 
\frac{  R^{1-\frac{m(n+p)}{p}\left (\frac{1}{p}-\frac{1}{q}  \right )} v \left (Q_R^{x_0} \right )^{\frac{1}{q}}}{ \omega \left (Q_R^{x_0}  \right )^{\frac{1}{p}-\frac{m}{p}\left (\frac{1}{p}-\frac{1}{q}  \right )}} \Bigg )^q\left \Vert  u \right \Vert_E^q 
$$
$$
+\frac{1}{4p}\Vert u \Vert_E^p -\dfrac{\mu}{\gamma}C\vert \Omega \vert^{{1}-\frac{\gamma}{N^{,}}}\Vert1+\omega^\frac{-p^{,}}{p}\Vert_{L_{1}(\Omega)} ^\frac{\gamma}{p^{,}}\Vert u\Vert^{\gamma} _{E}
$$
i.e. 
$$I(u)\geq \Bigg ( \frac{1}{4p} -\frac{1}{q}\Big (C_0 
\frac{ R^{1-\frac{m(n+p)}{p}\left (\frac{1}{p}-\frac{1}{q}  \right )} v \left (Q_R^{x_0} \right )^{\frac{1}{q}}}{ \omega \left (Q_R^{x_0}  \right )^{\frac{1}{p}-\frac{m}{p}\left (\frac{1}{p}-\frac{1}{q}  \right )}} \Big )^q \Vert u \Vert_E^{q-p} \Bigg )  \Vert u \Vert_E^p
$$
$$
+\Vert u\Vert^{\gamma} _{E}\left(\frac{1}{4p}\Vert u \Vert_E^{p-\gamma} -\dfrac{\mu}{\gamma}C\vert \Omega \vert^{{1}-\frac{\gamma}{N^{,}}}\Vert1+\omega^\frac{-p^{,}}{p}\Vert_{L_{1}(\Omega)} ^\frac{\gamma}{p^{,}}\right)+\frac{1}{2p}\Vert u \Vert_E^p,
$$

If we choose a sphere in $ E $  as  
$$
\frac{1}{q}\Big (C_0 
\frac{ R^{1-\frac{m(n+p)}{p}\left (\frac{1}{p}-\frac{1}{q}  \right )} v \left (Q_R^{x_0} \right )^{\frac{1}{q}}}{ \omega \left (Q_R^{x_0}  \right )^{\frac{1}{p}-\frac{m}{p}\left (\frac{1}{p}-\frac{1}{q}  \right )}}\Big )^q \Vert u \Vert_E^{q-p}=\frac{1}{4p}.  
$$
Precisely, if we set radius of the ball $B(0,\rho) \subset E $ with
\begin{equation}\label{HL}
\rho = \left( \left ( \frac{q}{4p} \right )^{\frac{1}{q}}  \frac{\omega(Q_R^{x_0})^{\frac{1}{p}-\frac{m}{p}\left (\frac{1}{p}-\frac{1}{q}  \right )}}{C_{0} R^{1-\frac{m(n+p)}{p}\left (\frac{1}{p}-\frac{1}{q}  \right )} v \left (Q_R^{x_0} \right )^{\frac{1}{q}}}  \right) ^{\frac{q}{q-p}
},
\end{equation}
and let $\mu$ to be sufficiently small such that $\mu \in (0,\Lambda)$ with
$$
\Lambda =  \left( \left ( \frac{q}{4p} \right )^{\frac{1}{q}}  \frac{\omega(Q_R^{x_0})^{\frac{1}{p}-\frac{m}{p}\left (\frac{1}{p}-\frac{1}{q}  \right )}}{C_{0} R^{1-\frac{m(n+p)}{p}\left (\frac{1}{p}-\frac{1}{q}  \right )} v \left (Q_R^{x_0} \right )^{\frac{1}{q}}}  \right) ^{\frac{q(2-p)}{q-p}} \vert \Omega \vert^{\frac{\gamma}{N}-{1}}\Vert1+\omega^\frac{-p^{,}}{p}\Vert_{L_{1}(\Omega)} ^\frac{-\gamma}{p^{,}} $$
as a result, we get the following estimate
\begin{equation}\label{PG}
I(u)\geq \frac{1}{2p} \left( \left ( \frac{q}{4p} \right )^{\frac{1}{q}} \frac{\omega(Q_R^{x_0})^{\frac{1}{p}-\frac{m}{p}\left (\frac{1}{p}-\frac{1}{q}  \right )}}{C_{0} R^{1-\frac{m(n+p)}{p}\left (\frac{1}{p}-\frac{1}{q}  \right )}  v \left (Q_R^{x_0} \right )^{\frac{1}{q}} }  \right) ^{\frac{pq}{q-p}}
\end{equation}
on sphere $S(0,\rho) \subset E$.

Now, it remains to find a point $ u_{0}\in E $ lied out of the ball $ B(0,R) $ in $ E $ where $ I(u_{0})<0 $ .
To show it, apply the fibering method: for a fixed $ u\in E $ and
sufficiently large $ t>1 $ it holds%
$$
I(tu)=\frac{t^p}{p}\Vert u \Vert_E^p -\frac{ t^q }{q}\int\limits_\Omega v u_+^q(z) \, dz- \frac{\mu t^{\gamma} }{\gamma}\int\limits_\Omega  u_{+}^{\gamma}(z) \, dz<0 .
$$
Applying MPT,
there exists a point $ \tilde u\in E $ with $ I(\tilde u)=c $ and $ I^{\prime }(\tilde u)=0.  $ Here
$$
c=  \inf \limits_{\beta (t) \in \Gamma} \sup\limits_{ t \in [0, 1]}  \, \, \left\{ I \left ( \beta
(t)\right )\right \},
$$
where the infimum is taken all over the curves $\Gamma =\left \{ \beta (t) \right \}$,
\[
\beta : [0,1]\rightarrow E,
\beta \in C^{1}[0,1;E],\quad
\text{and such that} \quad  \beta (0)=0,\beta (1)=u_0 .
\]
Therefore, $ I(\tilde u)>0,I ^{\prime }(\tilde u)=0. $ 

\bigskip
\textit{Existence of the second solution.}  we will find another positive solution $ u_1 $ such that $ u_1\ne u_0 $, using the variational calculus approach. Here the $u_0$ is a MPT solution.  Let $B(0, \rho)$ be the metric ball in $E$ with radius $ \rho $ defined in \eqref{HL} on the surface of which the inequality \eqref{PG} is fulfilled. Show that on each small neighborhood of zero, $B(0, \epsilon)\subset E$ there exist a point $ u \in  B(0, \epsilon)$ such that $I(u)<0.$ To show this, set the elements $\left \{ u_t= t \bar u: \, 0<t<1 \right \} $ lied in the ball $B(0, \rho)$ in $E, $ where $\bar u \in \partial B(0, \rho)$ is any fixed element on the surface of the ball. 
For the elements $u_t$
$$
I(u_t)=\frac{t^p}{p}\Vert \bar u \Vert_E^p -\frac{ t^q }{q}\int\limits_\Omega v \bar u_+^q (z) \, dz- \frac{ \mu t^\gamma }{\gamma}\int\limits_\Omega \bar u_+^\gamma (z) \, dz
$$
whence we get 
$$
I(u_t)=t^\gamma\frac{1}{p}\left ( t^{p-\gamma}\Vert \bar u \Vert_E^p -\frac{ \mu }{\gamma}\int\limits_\Omega \bar u_+^\gamma (z) \, dz\right )- \frac{ t^q }{q}\int\limits_\Omega \bar v u_+^q (z) \, dz \leq -\frac{ t^q  }{q}\int\limits_\Omega \bar v u_+^q (z) \, dz
$$
as $$0< t<\left ( \frac{1}{\gamma\Vert \bar u  \Vert ^p }\left ( \int_\Omega  \bar u_+^\gamma \, dz \right )\right )^{1/(p-\gamma)}.$$
Note that, $\Vert v_t \Vert =t\rho $ and $I(v_t)<-\frac{ t^q }{q}\int\limits_\Omega \bar v u_+^q (z) \, dz $ for  $ 0<t<1. $ 

Hence we have found a sequence of elements $v_t \to  0$ with $\Vert u_t\Vert=\rho$ and $I(u_t)<0$ as $ t \to 0.$ On sphere $S(0, \rho)\subset E$ we have the estimate  $$\inf _{S(0, \rho)} I(u)\geq \frac{1}{2p} \left( \left ( \frac{q}{4p} \right )^{\frac{1}{q}} \frac{\omega(Q_R^{x_0})^{\frac{1}{p}-\frac{m}{p}\left (\frac{1}{p}-\frac{1}{q}  \right )}}{C_{0} R^{1-\frac{m(n+p)}{p}\left (\frac{1}{p}-\frac{1}{q}  \right )}  v \left (Q_R^{x_0} \right )^{\frac{1}{q}} }  \right) ^{\frac{pq}{q-p}}.$$ As well as the functional $I(u)$ is bounded from below  $\Vert u \Vert <\rho: $
$$
I(u)\geq \frac{1}{p}\Vert u \Vert_E^p -\frac{1}{q}\Bigg (C_0 
\frac{  R^{1-\frac{m(n+p)}{p}\left (\frac{1}{p}-\frac{1}{q}  \right )} v \left (Q_R^{x_0} \right )^{\frac{1}{q}}}{ \omega \left (Q_R^{x_0}  \right )^{\frac{1}{p}-\frac{m}{p}\left (\frac{1}{p}-\frac{1}{q}  \right )}} \Bigg )^q\left \Vert  u \right \Vert_E^q 
$$
$$
-\dfrac{\mu}{\gamma}C\vert \Omega \vert^{{1}-\frac{\gamma}{N^{,}}}\Vert1+\omega^\frac{-p^{,}}{p}\Vert_{L_{1}(\Omega)} ^\frac{\gamma}{p^{,}}\Vert u\Vert^{\gamma} _{E}\geq
$$
$$
-\dfrac{\mu R^{\gamma}}{\gamma}C\vert \Omega \vert^{{1}-\frac{\gamma}{N^{,}}}\Vert1+\omega^\frac{-p^{,}}{p}\Vert_{L_{1}(\Omega)} ^\frac{\gamma}{p^{,}}.
$$

Let $$\inf_{B(0,\rho)} I(u)=m,$$ evidently $m<0.$ There exist a sequence of elements $\{ u_n \}$ such that $I(u_n)\to m $ as $n\to \infty. $ Then as proved above, $u_n$ is a bounded sequence. There exists a subsequence of $\{u_n\}$ with $u_n \rightharpoonup u_1$ in the weak topology of $E$ for some $u_1 \in E.$ Since the functional $I(u)$ is weakly lower semicontinuous,
$$\underline{\lim}_{n\to \infty} I(u_n)\geq I(u_1)\geq m, $$
i.e. $I(u_1)=m, $ therefore $u_1$ is a local minimum of the functional $I(u). $ Therefore, $u_1$ is a critical point of the functional $I(u).$ Also,  $I(u_1)<0, $ therefore, the solutions $u_0$ and $u_1$ are different 
since $I(u_1)<0<I(u_0).$ In the next subsection, we will show that the both solutions $u_0, u_1$ are positive.

\bigskip

\textit{Positivity of solutions.} 
Let $u(z)$ be a weak solution of the equation \eqref{eq1}.  To show that $%
u  $ is positive, take  $ \phi=u_- $  in $ \langle I^{\prime }(u), \phi  \rangle =0:$
\begin{equation*}
\begin{split}
	\left \langle I^\prime(u), u_{-} \right \rangle &=\int\limits_{\Omega}\vert \nabla_\omega u \vert^{p-2}\nabla_\omega u\ \nabla_\omega u_{-}  \, dz \\
	&- \int\limits_\Omega v(x) u_+^{q-1} u_{-} \, dz-\mu\int\limits_\Omega  u_{+} ^{\gamma-1} u_{-} \, dz =0 . 
\end{split}
\end{equation*}
The second and third integral is zero since $ v u_{+}(z)^{q-1} u_{-}(z)=0  $ and $  u_{+}(z)^{\gamma-1} u_{-}(z)=0  $ in $ \Omega. $ Set $D_+=\left \{ z\in \Omega: u(z)>0 \right \}$ and $D_{-}=\Omega\setminus D_{+}.$ Set $u(z)=u_{+}(z)-u_{-}(z),$ where $u_+=u \chi_{D_+} , \, u_{-}=u \chi_{D_{-}}$.

Then
$$
\int\limits_{\Omega}\vert \nabla_\omega u \vert^{p-2}\nabla_\omega u_{-}\ \nabla_\omega u_{-}  \, dz=0
$$
From this
applying the inequality of Sobolev type Lemma \ref{L2}, we get
$$
\int_{\Omega} v(x) u_-^q (z) \, dz =0,
$$
therefore $u_-(z) \equiv 0.$ 
Which proves positivity of the weak solutions $u(z).$

\hspace{10cm} $\Box .$

\bigskip

This completes the proof of Theorem \ref{t1}.

\bigskip

\bigskip

\bigskip

\textbf{Acknowledgement.} We thank the referees for their kind comments inspired us for achieving more good final version.

\end{document}